\documentclass[12pt]{article}
 \usepackage{a4, latexsym, amsfonts,amsmath}
 \newtheorem{theorem}{Theorem}

 \newtheorem{definition}{Definition}

 \newtheorem{open problem}{Open problem}[section]

 \newenvironment{proof}{\trivlist
      \item[\hskip\labelsep
      {\itshape Proof.}]\normalfont}
      {\hspace*{\fill}$\Box$\endtrivlist}

\begin{document}

\title{Positive rational number of the form $\varphi(km^{a})/\varphi(ln^{b})$\\}

\author{Hongjian Li \thanks{is with School of Mathematics, South China Normal University, Guangzhou 510631, China (email: 2523253922@qq.com).  }\,\,  Pingzhi Yuan\thanks{P. Yuan is with School of Mathematics, South China Normal University, Guangzhou 510631, China (email: yuanpz@scnu.edu.cn).}\,\,   Hairong Bai \thanks{ is with School of Mathematical Science,South China Normal University,	Guangzhou 510631, China (email:baihairong2007@163.com)}}

\date{}
\maketitle
 \edef \tmp {\the \catcode`@}
   \catcode`@=11
   \def \@thefnmark {}

\begin{abstract}

Let $k, l, a$ and $b$ be positive integers with $\max\{a, \, b\}\ge2$. In this paper, we show that every positive rational number can be written as the form $\varphi(km^{a})/\varphi(ln^{b})$, where $m, \, n\in\mathbb{N}$ if and only if $\gcd(a, \,b)=1$ or $(a, b, k, l)=(2,2, 1, 1)$. Moreover, if $\gcd(a, b)>1$, then the proper representation of such representation  is unique.

{\bf Keywords:} $\quad$  representation of positive rational numbers, Euler's totient function.

{\bf 2010 Mathematics Subject Classification}. Primary 11A25; Secondary 11D85.
 \end{abstract}

Let $\mathbb{N}$ be the set of all positive integers and  $\varphi$  Euler's totient function, which is defined as $\varphi(n)=\sharp\{r: r\in \mathbb{N}, 0<r\le n, \gcd(r, n)=1\}$, the number of integers in the set $1, 2, \ldots, n$ that are relatively prime to $n$. Let $$n=\prod_{i=1}^{s}p_{i}^{\alpha_{i}}, \quad p_1<p_2<\cdots<p_s \,\, \mbox{ are} \,\,  \mbox{ primes}, \alpha_i\in\mathbb{N}$$ be the standard factorization of a positive integer $n$, it is well-known that \begin{equation}\label{eq1}
\varphi(n)=\varphi(\prod_{i=1}^{s}p_{i}^{\alpha_{i}})=\prod_{i=1}^{s}(p_{i}-1)p_{i}^{\alpha_{i}-1}.\end{equation}
(See. e.g., \cite{IK90}, page 20). Sun \cite{Sun1} proposed many challenging conjectures on representations
of positive rational numbers. Recently, D. Krachun and Z. Sun \cite{Sun2} proved that: any positive rational number can be written as the form $\varphi(m^{2})/\varphi(n^{2})$,  where $m, \, n\in\mathbb{N}$.

For given positive integers $a, b, k$ and $l$ with $\max\{a, \, b\}\ge2$, in this paper, we consider a more general problem when every positive rational number can be written as the form $\varphi(km^{a})/\varphi(ln^{b})$,  where $m, \, n\in\mathbb{N}$.

Note that if $p$ is a prime with either $p|\gcd(m, n)$ or $\gcd(p, \, mn)=1$, then $\frac{\varphi(m)}{\varphi(n)}=\frac{\varphi(mp)}{\varphi(np)}$. From this we can easily derive that
\begin{equation}\label{eq2}
\frac{\varphi(km^ad_1^a)}{\varphi(ln^bd_2^b)}=\frac{\varphi(km^a)}{\varphi(ln^b)},\end{equation}
whenever $d_1$ and $d_2$ are positive integers with $d_1^a=d_2^b$, and for each prime $p|d_1d_2$, either $p|\gcd(km, \,ln)$ or $\gcd(p, klmn)=1$. Hence we have the following definition.

\begin{definition} Let $k, l, a$ and $b$ be positive integers with $\min\{a, \, b\}\ge2$ and $r$ a positive rational numbers. A representation of $r=\frac{\varphi(km^a)}{\varphi(ln^b)}$ is called a $proper$ representation if there are no positive integers $d_1>1$ and $d_2$ such that $m=m_1d_1$, $n=n_1d_2$, $d_1^a=d_2^b$, and for each prime $p|d_1d_2$, either $p|\gcd(km_1, \,ln_1)$ or $\gcd(p, klm_1n_1)=1$.\end{definition}

For example, when $(a, b, k, l)=(2, 2, 1, 1)$, then  a representation of $r=\frac{\varphi(m^2)}{\varphi(n^2)}$ is called a $proper$ representation if there are no positive integers $d>1, m_1, n_1$ such that $m=m_1d$, $n=n_1d$, and for each prime $p|d$, either $p|\gcd(m_1, \,n_1)$ or $\gcd(p, m_1n_1)=1$.

The main purpose of this note is to show  the following result.
\begin{theorem}\label{thm1}
Let $k, l, a$ and $b$ be positive integers with $\max\{a, \, b\}\ge2$. Then any positive rational number can be written as the form $\varphi(km^{a})/\varphi(ln^{b})$,  where $m, \, n\in\mathbb{N}$ if and only if $\gcd(a, \,b)=1$ or  $(a, \,b, k, l)=(2, \,2, 1, 1)$. Moreover, if $\gcd(a, b)>1$, then the proper representation of such representation of a  positive rational number is unique.\end{theorem}

{\bf Proof of Theorem 1:}

\begin{proof} For a positive rational number $r$ with $r\ne1$, let $$r=\prod_{i=1}^{s}p_{i}^{\alpha_{i}}, \quad p_{1}<p_{2}<\cdots<p_{s},\,\, \alpha_{i}\in \mathbb{Z}\backslash\{0\}$$ be the standard factorization of  $r$, where $p_{1}<p_{2}<\cdots<p_{s}$ are primes. Let $P(r)=p_s$ denote the maximal prime factor of $r$, $v_{p_s}(r)=\alpha_s$, the $p_s$-valuation of $r$, and we let $P(1)=1$.

We first prove that if $\gcd(a, b)=1$ or $(a, b)=(2, 2)$ and $(k, \, l)=(1, \,1)$, then any positive rational number can be written as the form $\varphi(km^{a})/\varphi(ln^{b})$,  where $m, \, n\in\mathbb{N}$. Recall that Krachun and Sun \cite{Sun2} have proved the case of $(a, b, \, k, \, l)=(2, 2, 1,\, 1)$. It  suffices to show the statement holds  for $a$ and $b$ with $\gcd(a, \, b)=1$.
  For any  integer $c$, it is well-known that there are positive integers $x$ and $y$ such that $ax-by=c$ since $\gcd(a, \, b)=1$.  Hence for each prime $p$ with $v_p(r)\ne0$ or $p|kl$, there are positive integers $x(p)$ and $y(p)$ such that
$$v_p(r)=v_p(k)+ax(p)-v_p(l)-by(p).$$
Let
$$m=\prod_{p|kl, \,\,\mbox{or} \,\,v_p(r)\ne0}p^{x(p)}, \quad n=\prod_{p|kl,\,\, \mbox{or}\, \,v_p(r)\ne0}p^{y(p)}.$$ Then it is easy to check that
$$r=\frac{\varphi(km^{a})}{\varphi(ln^{b})}$$
since for any prime $p$ with $v_p(r)\ne0$ or $p|kl$, we have $p|\gcd(m,n)$ and hence
$$v_p(r)=v_p(k)+av_p(m)-1-(v_p(l)+bv_p(n)-1)=v_p(k)+ax(p)-v_p(l)-by(p).$$
Therefore we have proved  that any positive rational number can be written as the form $\varphi(km^{a})/\varphi(ln^{b}), m, \, n\in\mathbb{N}$ when $\gcd(a, b)=1$ or $(a, b)=(2, 2)$ and $(k, \, l)=(1, \,1)$.

Next, we will show that if $\gcd(a, \, b)\ge2$ and $(a, b, k, l)\ne (2, \,2, 1, 1)$, then there exists a positive rational number $r$ such that $r$ cannot be written as the form $$\frac{\varphi(km^{a})}{\varphi(ln^{b})}, \quad  m, \, n\in\mathbb{N}.$$

If $\gcd(a, b)=\mu>2$,  take $t$ to be a positive integer with $t\equiv 1\pmod{ \mu}$ and $p$  a prime with $p>P(kl)$, we will show that $p^t$ cannot be written as the form $\frac{\varphi(km^{a})}{\varphi(ln^{b})}$, $ m, \, n\in\mathbb{N}$.

Suppose that there are positive integers $m$ and $n$ such that
$$p^t=\frac{\varphi(km^{a})}{\varphi(ln^{b})}.$$
Without loss of generality, we may assume that the above representation is  proper, so $P(klmn)=p$, and hence $P(mn)=p$ since $p>P(kl)$. Now we have
$$
t=\begin{cases}
 v_p\left(\varphi(p^{av_p(m)})/\varphi(p^{bv_p(n)})\right)=av_p(m)-bv_p(n)\equiv 0\pmod{\mu}, & \mbox{if}\,\,v_p(n)\neq 0, \\
 v_p\left(\varphi(p^{av_p(m)})\right)=av_p(m)-1, & \mbox{if}\,\, v_p(n)=0, \\
\end{cases}    $$
which implies that $t\equiv0\pmod{\mu}$ or $t\equiv-1\pmod{\mu}$, a contradiction to $t\equiv 1\pmod{ \mu}$ and $\mu\ge3$. Hence
 $p^{t}$ can not be written as $$\frac{\varphi(km^{a})}{\varphi(ln^{b})}, \quad \, m, \, n\in\mathbb{N}.$$

  For the case of  $\gcd(a, b)=\mu=2$ and $(a, b)\ne(2, 2)$, we have that $\max\{a, b\}>2$. We only consider the case where $a>2$ (the argument for the case of $b>2$ is similar). We have, for a prime $p$ with $p>P(kl)$,
 $$ v_p\left(\frac{\varphi(km^{a})}{\varphi(ln^{b})}\right)=av_p(m)-bv_p(n)\equiv0\pmod{2}\equiv0\pmod{2}, $$
or $av_p(m)-1\ge3$ when $v_p(n)=0$,  or $1-bv_p(n)<0$ when $v_p(m)=0$, so
$p$ can not be written as
$$\frac{\varphi(km^{a})}{\varphi(ln^{b})}, \quad m, \, n\in\mathbb{N}.$$

Now we consider the case where  $\gcd(a, b)=(2, 2)$ and $(k, l)\ne(1, 1)$, then $\max\{k, l\}\ge2$. Without loss of generality, we may assume that $p=P(kl)=P(k)$ (the argument for the case of $p=P(kl)=P(l)$ is similar). If $p|\gcd(k, l)$, then we have
$$v_p\left(\frac{\varphi(km^{a})}{\varphi(ln^{b})}\right)=v_p(k)-v_p(l)+2(v_p(m)-v_p(n))\equiv v_p(k)-v_p(l)\pmod{2},$$
so $p^{|v_p(k)-v_p(l)|+1}$ can not be written as
$\frac{\varphi(km^{a})}{\varphi(ln^{b})}, m, \, n\in\mathbb{N} $. If $p\not|l$ (the case where $p\not|k$ is similar, and we omit the detail), then we have
$$v_p\left(\frac{\varphi(km^{a})}{\varphi(ln^{b})}\right)=v_p(k)+2(v_p(m)-v_p(n))\equiv v_p(k)\pmod{2},$$
  or $v_p(k)-1+2v_p(m)\ge0$ when $v_p(n)=0$, so $p^{-v_p(k)-1}$ can not be written as
$\frac{\varphi(km^{a})}{\varphi(ln^{b})}, m, \, n\in\mathbb{N} $. This proves the first statement.

To   prove the last statement, we use a double induction on $P(kl)$. We first prove that the statement holds for $P(kl)=1$ by induction on $P(r)$, then we prove that the statement holds for any positive integers $k, l$ by induction on $P(kl)$.

 Let $d=\gcd(a, b)>1$ and $r$ a positive rational number.
Suppose that the representation of
\begin{equation}\label{eq21}
r=\frac{\varphi(km^a)}{\varphi(ln^b)},\quad \, m, \, n\in\mathbb{N}.\end{equation}
is proper. We will show that $m$ and $n$ are uniquely determined by $r$, $k$ and $l$, in other words, the proper representation (\ref{eq21}) of $r$ is unique. We prove this by a double induction on $P(kl)$. To begin with, we show that for the proper representation (\ref{eq21}) of $r$ is unique when $P(kl)=1$, i.e., $k=l=1$. In this case, we use induction on $P(r)$. For $P(r)=1$, let $1=\frac{\varphi(m^a)}{\varphi(n^b)}$ be a proper representation of $1$, and let $p=P(mn)$. If $p>1$,  then $p$ is a prime and we have
$$0=v_p(1)=av_p(m)-bv_p(n).$$
Hence $1=\frac{\varphi(m_1^a)}{\varphi(n_1^b)}$, where $m_1=m/p^{v_p(m)}$, $n_1=n/p^{v_p(n)}$,  and $av_p(m)=bv_p(n)$, which implies that the representation of $1=\frac{\varphi(m^a)}{\varphi(n^b)}$ is not proper, a contradiction. Hence $p=1$, and $1$ has only the proper representation of $1=\frac{\varphi(1^a)}{\varphi(1^b)}$. For a proper representation (\ref{eq21}) of $r$ with $P(r)>1$,  we claim that

\begin{equation}\label{eq22}
P(mn)=P(r).\end{equation}
Let $q= P(mn)$. Obviously,  we have $q\ge P(r)$. Note that
\begin{equation}\label{eq23} v_q\left(\frac{\varphi(m^a)}{\varphi(n^b)}\right)=\begin{cases}
 av_q(m)-bv_q(n), & \mbox{if}\, \, v_q(m)\ne0\,\, \mbox{and}\,\,  v_q(n)\ne0, \\
 av_q(m)-1, &  \mbox{if}\, \, v_q(n)=0, \\
 -bv_q(n)+1, & \mbox{if} \,\, v_q(m)=0.\\
\end{cases}  \end{equation}
If $q>P(r)$, then we have $0=v_q(r)=v_q\left(\frac{\varphi(m^a)}{\varphi(n^b)}\right)$, so $v_q(m)\ne0$ and $v_q(n)\ne0$. Hence $0=av_q(m)-bv_q(n)$, which implies that $av_q(m)=bv_q(n)$. Therefore, we have
$$r=\frac{\varphi(m^a)}{\varphi(n^b)}=\frac{\varphi(m_1^a)}{\varphi(n_1^b)},$$
where $m_1=m/q^{v_q(m)}$ and $n_1=n/q^{v_q(n)}$, i.e., the representation of $r=\frac{\varphi(m^a)}{\varphi(n^b)}$ is not proper by definition, a contradiction. Hence $P(mn)=P(r)$ when the representation (\ref{eq21}) of $r$ is proper.

For $P(r)=2$, by (\ref{eq22}), for any nonzero integer $c$ and any proper representation of
$$2^c=\frac{\varphi(m^a)}{\varphi(n^b)}, \quad  m, \, n\in\mathbb{N}$$
we have that $P(mn)=2$, so there are non-negative integers $x$ and $y$ such that $m=2^x$ and $n=2^y$. For simplicity, we assume that $c>0$ (the case where $c<0$ is similar). If $d=\gcd(a, b)|c$, then it follows from (\ref{eq23}) that $c=ax-by$, $x>0, y>0$. Let $(x, y)=(x_0, y_0)$, $x_0>0, \, y_0>0$ be the least positive integer solution of the linear equation
$$ax-by=c, \quad x>0, y>0.$$ It is well-known that all positive integer solution $(X, Y)$ of the equation $aX-bY=c$ are given by
\begin{equation}\label{l1}
X=x_0+\frac{b}{\gcd(a,b)}t, \quad Y=y_0+\frac{a}{\gcd(a,b)}t,\quad t\in\mathbb{N}\cup\{0\}.\end{equation} We claim that $x=x_0$ and $y=y_0$. Otherwise, there exists a positive integer $t$ such that
$$x=x_0+\frac{b}{\gcd(a,b)}t, \quad y=y_0+\frac{a}{\gcd(a,b)}t.$$ Hence we have
$$2^c=\frac{\varphi(2^{ax}d_1^a)}{\varphi(2^{by}d_2^b)}=\frac{\varphi(2^{ax_0})}{\varphi(2^{by_0})},\quad d_1^a=d_2^b=2^{\frac{abt}{\gcd(a, b)}},$$
which contradicts with the proper representation of $2^c$. Therefore  $x=x_0$ and $y=y_0$. Hence the proper representation of $2^c$ is unique. If $\gcd(a, b)\not|c$, then it follows from (\ref{eq23}) that $c=ax-1$ and $y=0$, so $x=(c+1)/a$, the proper representation of $2^c$ is also unique. This proves the result for $P(r)=2$.

Now let $q$ be an odd prime and assume that the result holds for $P(r)<q$. Let $r$ be a positive rational number with $P(r)=q$, and the representation (\ref{eq21}) of $r$ is proper. By (\ref{eq22}), we have $P(mn)=q$. If $v_q(n)=0$, then $v_q(r)>0$ and it follows from (\ref{eq23}) that $v_q(r)=av_q(m)-1$ and $v_q(n)=0$, so $v_q(m)=(v_q(r)+1)/a$. Let $r_0=r/q^{v_q(r)}(q-1)$. Then we have
$$r_0=\frac{\varphi(m_1^a)}{\varphi(n_1^b)},$$
where $m_1=m/q^{v_q(m)}$ and $n_1=n$. Note that the representation of $r_0=\frac{\varphi(m_1^a)}{\varphi(n_1^b)}$ is proper, so by the induction hypothesis, $m_1$ and $n_1$ are uniquely determined by $r_0$ since $P(r_0)<P(r)=q$. Therefore the proper representation of $r$ is unique.

If $v_q(m)=0$ and $v_q(n)>0$, then $v_q(r)<0$ and it follows from (\ref{eq23}) that $v_q(r)=-bv_q(n)+1$ and $v_q(m)=0$, so $v_q(n)=(-v_q(r)+1)/b$. Let $r_0=\frac{r(q-1)}{q^{v_q(r)}}$. Then we have
$$r_0=\frac{\varphi(m_1^a)}{\varphi(n_1^b)},$$
where $m_1=m$ and $n_1=n/q^{v_q(n)}$. Now the representation of $r_0=\frac{\varphi(m_1^a)}{\varphi(n_1^b)}$ is proper, so by the induction hypothesis  $m_1$ and $n_1$ is uniquely determined by $r_0$ since $P(r_0)<P(r)=q$. Therefore the proper representation of $r$ is unique.

 If $v_q(m)>0$ and $v_q(n)>0$,  then it follows from (\ref{eq23}) that $v_q(r)=av_q(m)-bv_q(n)$ since $v_q(m)>0$ and $v_q(n)>0$. By the similar argument as above, we have $(v_q(m), v_q(n))=(x, y)$ is the least positive integer solution of the linear equation
$ax-by=v_q(r)$. Let $r_0=r/q^{v_q(r)}$. Then we have
$$r_0=\frac{\varphi(m_1^a)}{\varphi(n_1^b)},$$
where $m_1=m/q^{v_q(m)}$ and $n_1=n/q^{v_q(n)}$. It is easy to verify that the representation of $r_0=\frac{\varphi(m_1^a)}{\varphi(n_1^b)}$ is proper.  Since $P(r_0)<P(r)=q$, by the induction hypothesis, $m_1$ and $n_1$ are uniquely determined by $r_0$. Therefore the proper representation of $r$ is unique.

In view of the above, we have proved that the representation (\ref{eq21}) of $r$ is unique. This completes the proof of the case $(k, \, l)=(1, 1)$.

Now we assume that the statement holds for $P(kl)<p$, where $p$ is a prime. That is, for any positive integers $k, l$ with $P(kl)<p$ and any proper representation of
$$r=\frac{\varphi(km^a)}{\varphi(ln^b)}, \quad  m, \, n\in\mathbb{N}$$  of a rational number number $r$, $m, n$ are uniquely determined by $r, k, l$.  We prove that the statement holds for $P(kl)=p$. For the first step, we let $r_0$ be a positive rational number with least $P(r_0)$ such that
$$r_0=\frac{\varphi(km^a)}{\varphi(ln^b)}, \quad  m, \, n\in\mathbb{N}$$
is a proper representation of $r_0$. We show that $m, n$ are uniquely determined by $r_0, k, l$. Obviously, $q=P(klmn)\ge P(r_0)$. If $q=P(klmn)>P(r_0)$, then we have
$$ 0= v_q\left(\frac{\varphi(km^a)}{\varphi(ln^b)}\right)=\begin{cases}
 av_q(m)+v_q(k)-bv_q(n)-v_q(l), & \mbox{if}\, \, v_q(km)\ne0\,\, \mbox{and}\,\,  v_q(ln)\ne0, \\
 av_q(m)+v_q(k)-1, &  \mbox{if}\, \, v_q(ln)=0, \\
 -bv_q(n)-v_q(l)+1, & \mbox{if} \,\, v_q(km)=0.\\
\end{cases}  $$
If $v_q(ln)=0$ (resp. $v_q(km)=0$), then $v_q(n)=0$ and $v_q(m)$ is uniquely determined by $k$ (resp. $v_q(m)=0$ and $v_q(n)$ is uniquely determined $l$). If $v_q(km)\ne0$ and $v_q(ln)\ne0$, then we have $ 0=av_q(m)+v_q(k)-bv_q(n)-v_q(l)$. By the same argument as before, we see that $(v_q(m), v_q(n))=(x, y)$ is the least non-negative (resp. positive) solution of the linear equation $ax-by=v_q(l)-v_q(k)$ when $v_q(k)>0$ and $v_q(l)>0$ (resp. when $v_q(k)=0$ or $v_q(l)=0$ ). Hence $v_q(m), v_q(n)$ are uniquely determined by $k, l$. We have
$$r_1=r=\frac{\varphi(k_1m_1^a)}{\varphi(l_1n_1^b)},$$
where $k_1=k/q^{v_q(k)}$, $l_1=l/q^{v_q(l)}$, $m_1=m/q^{v_q(m)}$, and  $n_1=n/q^{v_q(n)}$, and the above representation of $r_1$ is proper.  Further, if $q>p=P(kl)$, then $(k, l)=(k_1, l_1)$ and $av_q(m)=bv_q(n)$, so the representation of $r_0=\frac{\varphi(km^a)}{\varphi(ln^b)}, \quad  m, \, n\in\mathbb{N}$ is not proper. Hence $q=p=P(kl)$, and therefore $P(k_1l_1)<P(kl)=p$. By the induction hypothesis on $P(kl)$, we know that the representation of $r_0=\frac{\varphi(k_1m_1^a)}{\varphi(l_1n_1^b)}$ is unique, so $m_1, n_1$ are uniquely determined by $r_0, k_1, l_1$, and hence $m, n$ are uniquely determined by $r_0, k, l$.

If $q=P(klmn)=P(r_0)$, then we have
$$ v_q(r_0)= v_q\left(\frac{\varphi(km^a)}{\varphi(ln^b)}\right)=\begin{cases}
 av_q(m)+v_q(k)-bv_q(n)-v_q(l), & \mbox{if}\, \, v_q(km)\ne0\,\, \mbox{and}\,\,  v_q(ln)\ne0, \\
 av_q(m)+v_q(k)-1, &  \mbox{if}\, \, v_q(ln)=0, \\
 -bv_q(n)-v_q(l)+1, & \mbox{if} \,\, v_q(km)=0.\\
\end{cases} $$
By the same argument as in the case of $P(klmn)>P(r_0)$, we get  $v_q(m)$, $v_q(n)$ are uniquely determined by $r_0, k, l$  and we have
$$r_1=r_0/q^{v_q(r_0)}=\frac{\varphi(k_1m_1^a)}{\varphi(l_1n_1^b)},$$
where $k_1=k/q^{v_q(k)}$, $l_1=l/q^{v_q(l)}$, $m_1=m/q^{v_q(m)}$, and  $n_1=n/q^{v_q(n)}$, and the above representation of $r_1$ is proper. Further, we have $P(r_1)<P(r_0)$ and $P(k_1l_1)\le P(kl)=p$. If $P(k_1l_1)=P(kl)=p$, then $ (k, l)=(k_1, l_1)$,
$$r_1=\frac{\varphi(km_1^a)}{\varphi(ln_1^b)}$$
and $P(r_1)<P(r_0)$, which is impossible by the definition of $r_0$. If $P(k_1l_1)<P(kl)=p$, then by  the induction hypothesis on $P(kl)$, we know that the representation of $r_1=\frac{\varphi(k_1m_1^a)}{\varphi(l_1n_1^b)}$ is unique, so $m_1, n_1$ are uniquely determined by $r_1, k_1, l_1$, and hence $m, n$ are uniquely determined by $r_0, k, l$.

Next we assume that the statement holds for all positive rational number $r$ with $P(r_0)\le P(r)<q$, $q$ is an odd prime. We show that the statement holds for all positive rational number $r$ with $P(r)=q$. Let
$$r=\frac{\varphi(km^a)}{\varphi(ln^b)}, \quad  m, \, n\in\mathbb{N}$$
be a proper representation of $r$. Obviously, $q=P(klmn)\ge P(r)$. If $q=P(klmn)>P(r)$, then we have
\begin{equation}\label{eq25} 0= v_q\left(\frac{\varphi(km^a)}{\varphi(ln^b)}\right)=\begin{cases}
 av_q(m)+v_q(k)-bv_q(n)-v_q(l), & \mbox{if}\, \, v_q(km)\ne0\,\, \mbox{and}\,\,  v_q(ln)\ne0, \\
 av_q(m)+v_q(k)-1, &  \mbox{if}\, \, v_q(ln)=0, \\
 -bv_q(n)-v_q(l)+1, & \mbox{if} \,\, v_q(km)=0.\\
\end{cases}  \end{equation}
By the same argument as in the case of $r=r_0$, $v_q(m), v_q(n)$ are uniquely determined by $k, l$. And we have
$$r_1=r=\frac{\varphi(k_1m_1^a)}{\varphi(l_1n_1^b)},$$
where $k_1=k/q^{v_q(k)}$, $l_1=l/q^{v_q(l)}$, $m_1=m/q^{v_q(m)}$, and  $n_1=n/q^{v_q(n)}$, and the above representation of $r_1$ is proper. Similarly,  by the induction hypothesis on $P(kl)$, we know that the representation is unique, so $m_1, n_1$ are uniquely determined by $r, k_1, l_1$, and hence $m, n$ are uniquely determined by $r, k, l$.

The argument of the case where $P(klmn)=P(r)$ is the same as in the case of $r=r_0$. Now we have $v_q(m)$, $v_q(n)$ are uniquely determined by $r, k, l$  and
$$r_1=r/q^{v_q(r)}=\frac{\varphi(k_1m_1^a)}{\varphi(l_1n_1^b)},$$
where $k_1=k/q^{v_q(k)}$, $l_1=l/q^{v_q(l)}$, $m_1=m/q^{v_q(m)}$, and  $n_1=n/q^{v_q(n)}$, and the above representation of $r_1$ is proper. Further, we have $P(r_1)<P(r)$ and  $P(k_1l_1)\le P(kl)=p$. If $P(k_1l_1)=P(kl)=p$, then $ (k, l)=(k_1, l_1)$,
$$r_1=\frac{\varphi(km_1^a)}{\varphi(ln_1^b)}$$
and $P(r_1)<P(r_0)$, by the induction hypothesis on $P(r)$, we obtain that $m_1, n_1$ are uniquely determined by $r_1, k_1, l_1$, and hence $m, n$ are uniquely determined by $r_0, k, l$. If $P(k_1l_1)<P(kl)=p$, then by  the induction hypothesis on $P(kl)$, we know that the representation of $r_1=\frac{\varphi(k_1m_1^a)}{\varphi(l_1n_1^b)}$ is unique, so $m_1, n_1$ are uniquely determined by $r_1, k_1, l_1$, and hence $m, n$ are uniquely determined by $r, k, l$.

This completes the proof.

\end{proof}
{\bf Two examples:} (1) For $(a, b)=(2, 3)$ and $r=5/11$, we have
$$\frac{5}{11}=\frac{\varphi(5^{4}11^2)}{\varphi(5^{3}\cdot 11^3)}=\frac{\varphi(5^{2}11^2)}{4\varphi(11^3)}=\frac{\varphi(55^2)}{\varphi(22^3)}.$$
 Both $\frac{\varphi(275^2)}{\varphi(55^3)}$ and $\frac{\varphi(55^2)}{\varphi(22^3)}$ are proper representations of $\frac{5}{11}$.

(2) For $(a, b)=(2, 2)$ and $r=19/47$, we have
$$\frac{19}{47}=\frac{19\times46}{\varphi(47^2)}=\frac{19\varphi(23^2)}{11\varphi(47^2)}=\frac{5\varphi(19^2\times23^2)}{9\varphi(11^2\times47^2)}=
\frac{\varphi(3^2\times5^2\times19^2\times23^2)}{4\varphi(3^4\times11^2\times47^2)}$$
$$=\frac{\varphi(2^2\times3^2\times5^2\times19^2\times23^2)}{\varphi(2^4\times3^4\times11^2\times47^2)}=\frac{\varphi(13110^2)}{\varphi(18612^2)}=\frac{\varphi(39330^2)}{\varphi(55836^2)}.$$
Here $\frac{\varphi(13110^2)}{\varphi(18612^2)}$ is the proper representation of $r=19/47$, while the representation $\frac{\varphi(39330^2)}{\varphi(55836^2)}$ in \cite{Sun2} is not proper.

{\bf Remarks:} (1) It follows from the proof of Theorem 1 that if $a, b, m, n$ are positive integers such that $\min\{a, b\}>1$ and $\varphi(m^a)=\varphi(n^b)$, then $m^a=n^b$.

(2)  For $(a, b)=(2, 2)$, by the main theorem of \cite{Sun2} and Theorem 1, we have that any positive rational number can be written as the form $\varphi(m^{2})/\varphi(n^{2})$,  where $m, \, n\in\mathbb{N}$. Moreover, the proper representation of any positive rational number of the form $\varphi(m^2)/\varphi(n^2)$  is unique.

(3) Let $p$ be a prime and $a, b, c$ positive integers with $\gcd(a, \, b)=1$ and $c=au-1, u\in\mathbb{N}$. Let $\frac{1}{p-1}=\frac{\varphi(m^a)}{\varphi(n^b)}$ be a proper representation of $\frac{1}{p-1}$ and $x, y$ be the least positive integer solution of the linear equation $ax-by=au-1$. Then
$$p^c=\frac{\varphi(p^{xa})}{\varphi(p^{yb})}=\frac{\varphi(p^{ua}m^a)}{\varphi(n^{b})}.$$
And both representations above of $p^c$ are proper. The first equality of the above equation  does not hold when $\gcd(a, b)>1$ since the related linear equation $ax-by=au-1$ has no integer solution $(x, \, y)$.

{\bf Acknowledgments:}
This work is  supported by the National
Natural Science Foundation of China (grants no.12171163).


\begin{thebibliography}{99}
\bibitem{IK90} Ireland, K., Rosen, M.(1990).  A Classical Introduction to Modern Number Theory, 2nd ed. Graduate Texts in Mathematics, Vol. 84. New York: Springer.
\bibitem{Sun1} Sun, Z.-W.(2017). Conjectures on representations involving primes. In: Nathanson. M. ed. Combinatorial and Additive Number Theory II. Springer proceedings in mathematics and Statistics, Vol. 220. Cham: Springer, pp. 279-310.
\bibitem{Sun2} Krachun, D., Sun, Z.-W.(2020). Each positive rational number has the form $\varphi(m^{2})/\varphi(n^{2})$, Amer. Math. Monthly (to appear).
\end{thebibliography}
 \end{document}